\documentclass[11pt, reqno]{amsart}

\usepackage{amsmath, amsfonts, amssymb, amsbsy, bigstrut, graphicx, enumerate,  upref, longtable, comment, comment}

\usepackage{pstricks}

\usepackage[breaklinks]{hyperref}

\usepackage[numbers, sort&compress]{natbib}

\newcommand{\ep}{\epsilon}

\newtheorem{thm}{Theorem}[section]

\newtheorem{problem}[thm]{Problem}
\theoremstyle{definition}

\newcommand{\ee}{\mathbb{E}}

\newcommand{\ma}{\mathcal{A}}

\newcommand{\rr}{\mathbb{R}}
\newcommand{\smallavg}[1]{\langle #1 \rangle}

\newcommand{\tr}{\operatorname{Tr}}

\newcommand{\zz}{\mathbb{Z}}

\newcommand{\fpar}[2]{\frac{\partial #1}{\partial #2}}

\numberwithin{equation}{section}

\renewcommand{\Re}{\operatorname{Re}}

\begin{document}
\title{Yang--Mills for probabilists}
\author{Sourav Chatterjee}
\address{\newline Department of Statistics \newline Stanford University\newline Sequoia Hall, 390 Serra Mall \newline Stanford, CA 94305\newline \newline \textup{\tt souravc@stanford.edu}}
\thanks{Research partially supported by NSF grant DMS-1608249}

\dedicatory{Dedicated to friend and teacher Raghu Varadhan on the occasion of his $75^{th}$ birthday.}

\begin{abstract}
The rigorous construction of quantum Yang--Mills theories, especially in dimension four, is one of the central open problems of mathematical physics. Construction of Euclidean Yang--Mills theories is the first step towards this goal. This article presents a formulation of some of the core aspects this problem as problems in probability theory. The presentation begins with an introduction to the basic setup of Euclidean Yang--Mills theories and lattice gauge theories. This is followed by a discussion of what is meant by a continuum limit of lattice gauge theories from the point of view of theoretical physicists. Some of the main  issues are then posed as problems in probability. The article ends with a brief review of the mathematical literature.  \end{abstract}

\maketitle

\tableofcontents

\section{Introduction}\label{backsec}
Four-dimensional quantum Yang--Mills theories are the building blocks of the Standard Model of quantum mechanics. The Standard Model encapsulates the sum total of  all that is currently known about the basic particles of nature (see~\cite{griffiths08} for more details about the physics). Unfortunately, in spite of their incredible importance, quantum Yang--Mills  theories have no rigorous mathematical foundation. In fact, the mathematical foundation is so shaky that we do not even know for certain the right spaces on which these theories should be defined, or the right observables to look at.

Quantum Yang--Mills theories are defined in Minkowski spacetime. Euclidean Yang--Mills theories are `Wick-rotated' quantum Yang--Mills theories that are defined in Euclidean spacetime. These should, at least in principle, be easier to understand and analyze than their Minkowski counterparts. Although these theories are not rigorously defined either, they are formally probability measures on spaces of connections on certain principal bundles (more details to follow). Moreover, they have lattice analogs, known as lattice gauge theories or lattice Yang--Mills theories, that are rigorously defined probabilistic models. Therefore, the construction of Euclidean Yang--Mills theories, when viewed as a problem of taking scaling limits of lattice gauge theories, reveals itself as a problem in probability theory.

The problem of rigorously constructing Euclidean Yang--Mills theories, and then  extending the definition to Minkowski spacetime via Wick rotation, is the problem of Yang--Mills existence, posed as a `millennium prize problem' by the Clay Institute~\cite{jaffewitten}.

Quantum Yang--Mills theories are certain kinds of quantum field theories. A  standard approach to the rigorous construction of quantum field theories is via  the program of constructive quantum field theory, as outlined in the classic monograph of \citet{gj87}. One of the important objectives of constructive quantum field theory is to define Euclidean quantum field theories as probability measures on appropriate spaces of generalized functions, and then show that these probability measures satisfy certain axioms (the Wightman axioms, or the Osterwalder--Schrader axioms), which would then imply that the theory can be `quantized' to obtain the desired quantum field theories in Minkowski spacetime. However, as noted in the monograph of~\citet{seiler82}, there is a fundamental problem in following this path for Yang--Mills theories: the key observables in these theories do not take values at points, but at curves. Thus, it is not clear how to describe these theories as probability measures on spaces of generalized functions on manifolds --- they should, rather, be probability measures on spaces of generalized functions on spaces of curves. These considerations led \citet{seiler82} to pose the  problem as a problem of constructing an appropriate random function on a suitable space of closed curves. This is the route that we will adopt in this section, partly because it gives the most straightforward way of stating the problem.

That said, however, we still need a bit of preparation.  After presenting quick introductions to Euclidean Yang--Mills theories, lattice gauge theories and Wilson loops in the first three sections, the physicist's definition of a continuum limit of lattice gauge theories is discussed and well-defined mathematical problems are formulated. The last section contains a brief review of the mathematical literature.


\section{Euclidean Yang--Mills theories}\label{ymintro}
An Euclidean Yang--Mills theory involves a dimension $n$, and a `gauge group' $G$, usually a compact Lie group. For simplicity, let us assume that $G$ is a closed subgroup of the group of unitary matrices of some order $N$. Examples are $G=U(1)$ for quantum electrodynamics, $G=SU(3)$ for quantum chromodynamics, and $G = SU(3)\times SU(2)\times U(1)$ for the Standard Model, with dimension $n=4$ in each case. 

The Lie algebra $\mathfrak{g}$ of the Lie group $G$ is a subspace of the vector space of all $N\times N$ skew-Hermitian matrices. A $G$ connection form on $\rr^n$ is a smooth map from $\rr^n$ into $\mathfrak{g}^n$. If $A$ is a $G$ connection form, its value $A(x)$ at a point $x$ is an $n$-tuple $(A_1(x),\ldots, A_n(x))$ of skew-Hermitian matrices. In the language of differential forms, 
\[
A = \sum_{j=1}^n A_j dx_j\,.
\]
The curvature form $F$ of a connection form $A$ is the $\mathfrak{g}$-valued $2$-form
\[
F = dA + A \wedge A\,.
\]
In coordinates, this means that at each point $x$, $F(x)$ is an $n\times n$ array of skew-Hermitian matrices of order $N$, whose $(j,k)^{\mathrm{th}}$ entry is the matrix
\[
F_{jk}(x) = \fpar{A_k}{x_j} - \fpar{A_j}{x_k} + [A_j(x), A_k(x)]\,,
\]
where $[B, C] = BC - CB$ denotes the commutator of $B$ and $C$.

Let $\ma$ be the space of all $G$ connection forms on $\rr^n$. The Yang--Mills action on this space is the function
\[
S_{\mathrm{YM}}(A) := -\int_{\rr^n} \tr(F\wedge *F) \,,
\]
where $F$ is the curvature form of $A$ and $*$ denotes the Hodge star operator, assuming that this integral is finite. Explicitly, this is
\begin{equation}\label{ymaction}
S_{\mathrm{YM}}(A) = -\int_{\rr^n}\sum_{j,k=1}^n \tr(F_{jk}(x)^2)\,dx \,.
\end{equation}
The Euclidean Yang--Mills theory with gauge group $G$ on $\rr^n$ is formally described as the probability measure
\[
d\mu(A) = \frac{1}{Z}\exp\biggl(-\frac{1}{4g^2}S_{\mathrm{YM}}(A)\biggr) dA\,,
\]
where $A$ belongs to the space $\ma$ of all $U(N)$ connection forms, $S_{\mathrm{YM}}$ is the Yang--Mills functional defined above,  
\[
dA = \prod_{j=1}^n \prod_{x\in \rr^n} d(A_j(x))
\]
is infinite-dimensional Lebesgue measure on $\ma$, $g$ is a parameter called the coupling strength, and $Z$ is the normalizing constant that makes this a probability measure. 

The above description of Euclidean Yang--Mills theory with gauge group $G$ is not directly mathematically meaningful because any simple-minded way of defining infinite-dimensional Lebesgue measure on $\ma$ would yield $Z=\infty$ and make it impossible to define $\mu$ as a probability measure. While it has been possible to circumvent this problem in roundabout ways and give rigorous meanings to similar descriptions of Brownian motion and various quantum field theories in dimensions two and three, Euclidean Yang--Mills theories have so far remained largely intractable. 

\section{Lattice gauge theories}\label{lgtintro}
In 1974, \citet{wilson74} proposed a discretization of Euclidean Yang--Mills theories. These are now known as lattice gauge theories or lattice Yang--Mills  theories. Let $G$ be as in the previous section. The lattice gauge theory with gauge group $G$ on a finite set $\Lambda \subseteq\zz^n$ is defined as follows. Suppose that for any two adjacent  
vertices $x,y\in \Lambda$, we have a unitary matrix $U(x,y)\in G$, with the constraint that $U(y,x)=U(x,y)^{-1}$. Let us call any such assignment of matrices to edges a `configuration'. Let $G(\Lambda)$ denote the set of all configurations. A square bounded by four edges is called a plaquette. Let $P(\Lambda)$ denote the set of all plaquettes in $\Lambda$. For a plaquette $p\in P(\Lambda)$ with vertices $x_1,x_2,x_3,x_4$ in anti-clockwise order, and a configuration $U\in G(\Lambda)$, define
\begin{equation}\label{updef}
U_p := U(x_1,x_2)U(x_2,x_3)U(x_3,x_4)U(x_4,x_1).
\end{equation}
The Wilson action of $U$ is defined as 
\begin{align}\label{hdef}
S_\Lambda(U) := \sum_{p\in P(\Lambda)} \Re(\tr(I-U_p)),
\end{align}
where $I$ is the identity matrix of order $N$. 
Let $\sigma_{\Lambda}$ be the product Haar measure on $G(\Lambda)$. Given $\beta >0$, let $\mu_{\Lambda, \beta}$ be the probability measure on $G(\Lambda)$ defined as
\[
d\mu_{\Lambda, \beta}(U) := \frac{1}{Z}e^{-\beta S_\Lambda(U)} d\sigma_{\Lambda}(U)\, ,
\]
where  $Z$ is the normalizing constant. This probability measure is called the lattice gauge theory on $\Lambda$ for the gauge group $G$, with inverse coupling strength $\beta$. 

Often, it is convenient to work with an infinite volume limit of the theory, that is, a weak limit of the above probability measures as $\Lambda \uparrow \zz^n$. The infinite volume limit may or may not be unique. Indeed, the  uniqueness (or non-uniqueness) is in general unknown for lattice gauge theories in dimensions higher than two when $\beta$ is large.

Lattice gauge theories in several variations are a huge computational engine for numerical approximations for the Standard Model. They are used to make very accurate predictions of quantities like the masses of hadrons~\cite{bazavovetal}.  For readers who want to learn more about the physics of lattice gauge theories, two textbook references are \cite{montvaymunster97} and \cite{gattringerlang10}. There are also two rather extensive scholarpedia articles on lattice gauge theories and lattice quantum field theory, respectively.

The passage from a lattice gauge theory to an Euclidean Yang--Mills theory is heuristically justified as follows. First,  discretize the space $\rr^n$ as the scaled lattice $\ep \zz^n$ for some small $\ep$. Next, take a $G$ connection form 
\[
A= \sum_{j=1}^n A_j dx_j.
\]
Let $e_1,\ldots, e_n$ denote the standard basis vectors of $\rr^n$. For a directed  edge $(x, x+\ep e_j)$ of $\ep\zz^n$, define
\[
U(x,x+\ep e_j) := e^{\ep A_j(x)},
\]
and let $U(x+\ep e_j, x) := U(x,x+\ep e_j)^{-1}$. 
This defines a configuration of unitary matrices assigned to directed edges of $\ep\zz^n$. For a plaquette $p$ in $\ep \zz^n$,  let $U_p$ be defined as in \eqref{updef}. Then a formal calculation using the Baker--Campbell--Hausdorff formula for products of matrix exponentials shows that when $\ep$ is small, 
\begin{align*}
\sum_p \Re(\tr(I-U_p))&\approx \frac{\ep^{4-n}}{4}S_{\mathrm{YM}} (A),
\end{align*}
where $S_{\mathrm{YM}}$ is the Yang--Mills action defined in \eqref{ymaction}. The calculation goes as follows. Take any $x\in \ep \zz^n$ and any $1\le j<k\le n$, and let 
\[
x_1 = x, \ x_2=x+\ep e_j, \ x_3 = x+\ep e_j + \ep e_k, \ x_4 = x+\ep e_k.
\]
Let $p$ be the plaquette formed by the vertices $x_1,x_2,x_3,x_4$. Let $U_p$ be defined as in \eqref{updef}. Then 
\begin{align*}
U_p &=  e^{\ep A_j(x_1)} e^{\ep A_k(x_2)} e^{-\ep A_j(x_4)} e^{-\ep A_k(x_1)}.
\end{align*}
Recall the Baker--Campbell--Hausdorff formula for products of matrix exponentials:
\begin{align*}
e^B e^C &= \exp\biggl(B + C + \frac{1}{2}[B,C] + \text{higher commutators}\biggr)\,.
\end{align*}
Iterating this gives, for any $m$ and any $B_1,\ldots,B_m$,
\begin{align*}
e^{B_1}\cdots e^{B_m} &= \exp\biggl(\sum_{a=1}^m B_a + \frac{1}{2}\sum_{1\le a<b\le m}[B_a,B_b] + \text{higher commutators}\biggr).
\end{align*}
Recall that the eigenvalues of a skew-Hermitian matrix are all purely imaginary, and that the commutator of two skew-Hermitian matrices is skew-Hermitian. Consequently, the term within the exponential on the right side of the above display is skew-Hermitian and therefore has a purely imaginary trace. This implies that if $N$ is the order of the matrices, if the entries of $B_1,\ldots,B_m$ are of order $\ep$ and if the entries of $B_1+\cdots+B_m$ are of order $\ep^2$, then
\begin{align*}
\Re(\tr(I-e^{B_1}\cdots e^{B_m})) &= -\frac{1}{2}\tr\biggl[\biggl(\sum_{a=1}^m B_a + \frac{1}{2}\sum_{1\le a<b\le m}[B_a,B_b]\biggr)^2\biggr]\\
&\qquad  + O(\ep^5),
\end{align*}
where the real part of the trace was replaced by the trace on the right because the square of a skew-Hermitian matrix has real eigenvalues. Writing
\[
A_k(x_2) = A_k(x+\ep e_j) = A_k(x)+\ep \fpar{A_k}{x_j} + O(\ep^2)
\]
and using a similar Taylor expansion for $A_j(x_4)$, we get
\begin{align*}
A_j(x_1)+ A_k(x_2)- A_j(x_4) -A_k(x_1) = \ep\biggl(\fpar{A_k}{x_j} - \fpar{A_j}{x_k}\biggr) + O(\ep^2).
\end{align*}
Combining the above observations gives
\begin{align*}
\Re(\tr(I-U_p)) &= -\frac{1}{2}\ep^4 \tr\biggl[\biggl(\fpar{A_k}{x_j} - \fpar{A_j}{x_k} + [A_j(x), A_k(x)]\biggr)^2\biggr] + O(\ep^5)\\
&= -\frac{1}{2}\ep^4\tr(F_{jk}(x)^2)+O(\ep^5).
\end{align*}
This gives the formal approximation
\begin{align*}
S(U)  &= \sum_p \Re(\tr(I-U_p))\\
&\approx -\frac{1}{4}\sum_{x\in \ep \zz^n}\sum_{j,k=1}^n\ep^4 \tr(F_{jk}(x)^2)\\
&\approx -\frac{\ep^{4-n}}{4}\int_{\rr^n}\sum_{j,k=1}^n \tr(F_{jk}(x)^2)\, dx = \frac{\ep^{4-n}}{4}S_{\mathrm{YM}} (A).
\end{align*}
The above heuristic was used by Wilson to justify the approximation of Euclidean Yang--Mills theory by lattice gauge theory, scaling the inverse coupling strength $\beta$ like $\ep^{4-n}$ as the lattice spacing $\ep\to 0$. The most important dimension is $n = 4$, because spacetime is four-dimensional. In the above formulation, $\beta$ does not scale with $\ep$ at all when $n=4$. Currently, however, the general belief in the physics community is that $\beta$ should scale like $\log (1/\ep)$ in dimension four, although there are doubts about this belief and the question remains an open mathematical problem. 

The problem of constructing continuum limits of lattice gauge theories, from the point of view of theoretical physicists, is discussed in greater detail in the following sections. Most of this is `common knowledge' in the theoretical physics community, but not formalized in the sense of rigorous mathematics. I do not have references, but I have found \cite{kogut79, seiler82, kogut83} helpful. In particular, \cite{seiler82} proposes a formulation of continuum limits in terms of Wilson loop expectations, from which I borrow. 

\section{Wilson loop variables and quark confinement}\label{confinesec}
Any physical theory should have observables of interest. For Yang--Mills theories, the most important observables are Wilson loop variables. These are defined as follows. Suppose that we have an Euclidean Yang--Mills theory on $\rr^n$ with gauge group $G$, as defined in Section \ref{ymintro}. Given a piecewise smooth closed path $\gamma$ in $\rr^n$ and a $G$ connection $A$, the Wilson loop variable for $\gamma$ is defined as 
\[
W_\gamma:= \tr\biggl(\mathcal{P} \exp\biggl( \int_\gamma \sum_{j=1}^n A_j dx_j\biggr)\biggr),
\]
where $\mathcal{P}$ is the path-ordering operator.  In differential geometric terminology, the term inside the trace in the above display is the holonomy of $A$ along the closed path $\gamma$. Alternatively, it is the parallel transport of the identity matrix along $\gamma$ by the connection $A$. If the reader is unfamiliar with these concepts, there is nothing to worry. A simple definition of Wilson loop variables for lattice gauge theories is given below.

The physical importance of Wilson loop variables stems in part from their connection with the static quark potential. It was argued by \citet{wilson74} that the potential between a static quark and antiquark separated by distance $R$ is given by the formula
\[
V(R) = -\lim_{T\to\infty} \frac{1}{T}\log \smallavg{W_{\gamma_{R,T}}},
\]
where $\gamma_{R,T}$ is the boundary of a rectangle of length $T$ and breadth $R$, and $\smallavg{\cdot}$ denotes expectation with respect to a suitable Yang--Mills theory.  If $V(R)$ grows to infinity as $R\to\infty$, the quark-antiquark pair cannot separate beyond a fixed distance. This is the phenomenon of quark confinement, observed in experiments but currently lacking a satisfactory theoretical understanding (much less proof) due to the uncertainty about the existence of a continuum limit (although extensive numerical work in the lattice community points to a positive answer). In fact, it is believed that $V(R)$ grows like a multiple of $R$ for non-Abelian Yang--Mills theories in dimension four. This is known as Wilson's area law. If the area law holds, then the quantity 
\[
\lim_{R\to\infty}\frac{V(R)}{R}
\]
has physical significance. It is called the `string tension' of the continuum theory, and represents the energy density per unit length in the theory.


For lattice gauge theories, the definition of a Wilson loop variable is very simple. Suppose that we have a lattice gauge theory on  $\Lambda \subseteq \zz^n$ with gauge group $G$, as in Section \ref{lgtintro}.  A loop in $\zz^n$ is simply a directed path in the lattice which ends where it started. Given a loop $\gamma$ with directed edges $e_1,\ldots, e_m$, the Wilson loop variable $W_\gamma$ is defined as
\[
W_\gamma := \tr(U(e_1)U(e_2)\cdots U(e_m)). 
\]
The rationale for this definition is as follows. Let $A$ be a smooth $G$ connection on $\rr^n$. Take some small $\ep$ and define a configuration of group elements assigned to directed edges of $\ep\zz^n$ using the connection $A$, as in Section~\ref{lgtintro}. Let $\gamma$ be a smooth closed path in $\rr^n$ and let $\gamma_\ep$ be a loop in $\ep\zz^n$ that approximates this path. Then, as $\ep\to 0$, the discrete Wilson loop variable $W_{\gamma_\ep}$ approaches the continuous Wilson loop variable $W_\gamma$. 

Wilson's original motivation for investigating lattice gauge theories was to gain a theoretical proof of quark confinement. Since the expected values of Wilson loop variables are mathematically well-defined for lattice gauge theories, one can hope to give a rigorous proof of  the area law in the discrete setting. In fact, the following upper bound suffices. 
\begin{problem}[Area law]\label{pblmarea}
Take any compact non-Abelian Lie group $G\subseteq U(N)$ for some $N\ge 2$ and consider any infinite volume limit of four-dimensional lattice gauge theory with gauge group $G$ at inverse coupling strength $\beta$. Let $\gamma_{R,T}$ be a rectangular loop of breadth $R$ and length $T$ in the lattice. Prove that
\[
|\smallavg{W_{\gamma_{R,T}}}| \le C(\beta)e^{-c(\beta)RT },
\]
where $C(\beta)$ and $c(\beta)$ are positive constants that depend only on the inverse coupling strength $\beta$ and the gauge group.
\end{problem}
Soon after the appearance of Wilson's paper, physicists realized that the area law holds for any lattice gauge theory at sufficiently small $\beta$. A rigorous proof was given by \citet{os78}. However, this also implied that the area law at small $\beta$ cannot be evidence for quark confinement, because there are certain Yang--Mills theories that should not be confining quarks. An example is four-dimensional $U(1)$ Yang--Mills theory, which is the theory of electromagnetism. It is a fact of nature that there are no confined quarks in electromagnetism.

This apparent paradox was resolved by \citet{guth80}, who showed that four-dimensional $U(1)$ lattice gauge theory fails to satisfy the area law at large $\beta$. A fully rigorous proof of Guth's theorem was given by \citet{frohlichspencer82}. This result suggested that to prove quark confinement in a Yang--Mills theory, the area law has to be proved for the corresponding lattice theory at arbitrarily large $\beta$. For non-Abelian theories, this is currently known only in dimension two, where it is not very hard to prove. For $U(1)$ theory, there is a remarkable result of \citet{gopfertmack82}, who proved the area law at arbitrary $\beta$ in dimension three. The solution of Problem~\ref{pblmarea} for four-dimensional non-Abelian lattice gauge theories at large $\beta$ remains elusive.

\section{The problem of defining the continuum limit}\label{contsec}
Take a lattice gauge theory as in Section \ref{lgtintro}, and consider an infinite volume limit of this theory on $\zz^n$ obtained by taking a weak limit of the theories on finite cubes. Let $\gamma_1$ and $\gamma_2$ be two Wilson loops of fixed length, such as two plaquettes. The correlation between $W_{\gamma_1}$ and $W_{\gamma_2}$ is defined as the quantity
\begin{equation}\label{correq}
\smallavg{W_{\gamma_1}W_{\gamma_2}} - \smallavg{W_{\gamma_1}}\smallavg{W_{\gamma_2}}. 
\end{equation}
Let $d(\gamma_1,\gamma_2)$ denote the Euclidean distance between the two loops. If the logarithm of the above correlation behaves like $-d(\gamma_1,\gamma_2)/\xi$ for some $\xi>0$ as $d(\gamma_1,\gamma_2)\to\infty$, then the number $\xi$ is called the correlation length of the model. We have to take the logarithm because there may be polynomial correction terms to the exponential decay in the actual correlation~\cite{paesleme78, fariadaveigaetal04}.

Physicists say that the model has a continuum limit if there is a critical point $\beta_c\in [0,\infty]$ such that as $\beta \to\beta_c$, the correlation length tends to infinity. The reason for saying this is that if such a critical  point exists, then it is possible to define the model on the scaled lattice $\ep\zz^n$ instead of $\zz^n$, and send $\ep\to0$ in an appropriate manner as $\beta\to\beta_c$ such that the correlation length tends to a finite nonzero limit.  In other words, it is possible to define correlations in the continuum. 

It is believed that in dimension four (which, as stated earlier, is the dimension of greatest physical significance since spacetime is four-dimensional), many of the non-Abelian lattice models of interest have $\beta_c=\infty$. That is, one needs to take $\beta\to\infty$ while sending the lattice spacing $\ep$ to zero, to obtain a nontrivial correlation function in the limit. 

The following is a possible formulation of the above discussion as a concrete mathematical problem.  
\begin{problem}[Mass gap]\label{pblm0}
Take any compact non-Abelian Lie group $G\subseteq U(N)$ for some $N\ge 2$ and consider any infinite volume limit of four-dimensional lattice gauge theory with gauge group $G$ at inverse coupling strength $\beta$.  For each $x\in \rr^4$, let $p_x$ be the plaquette that is closest to $x$ (breaking ties by some arbitrary rule). Let $f_\beta(x)$ denote the correlation between $W_{p_0}$ and $W_{p_{x}}$, as defined in \eqref{correq}.  Show that for any $\beta>0$, there exists some $\xi(\beta)\in (0,\infty)$ such that 
\[
\lim_{|x|\to \infty} \frac{\log f_\beta(x)}{|x|} = -\frac{1}{\xi(\beta)}.
\]
Moreover, prove that 
\[
\lim_{\beta\to\infty} \xi(\beta)=\infty.
\]
\end{problem}

The correlation length $\xi$ has a physical meaning. Any lattice gauge theory contains information of an associated class of elementary particles called `glueballs' or `gluon-balls'. The existence of glueballs is considered to be one of the most important predictions of the Standard Model, but has not yet been experimentally verified. The number $\xi$ represents the reciprocal of the mass of the lightest glueball in the theory. 

One approach to the construction of continuum limits of lattice gauge theories is via Wilson loops. This is the approach that is advocated by~\citet[Chapter 8]{seiler82}, from which we draw inspiration. While~\citet{seiler82} gives a detailed description of the desired properties of the continuum limit that would presumably facilitate the quantization of the theory, we will restrict attention to the most basic question that needs to be solved before making any further progress.

Let $\beta_c$ be as above. The problem of constructing a continuum limit at this critical point in terms of Wilson loop expectations can be stated as follows. As $\beta \to \beta_c$,  one would like to show that  the lattice spacing  $\ep$ can be taken to $0$ in such a way that if $\gamma_\ep$ is any sequence of lattice loops converging to a  loop $\gamma$ in $\rr^n$, then $\smallavg{W_{\gamma_\ep}}$ converges to a nontrivial limit after some appropriate renormalization. 

Since $\infty$ is believed to be a critical point of compact non-Abelian lattice gauge theories in dimension four, one way to formulate the above question for the simple case of rectangular loops is the following.
\begin{problem}[Continuum limit]\label{pblm1}
Take any compact non-Abelian Lie group $G\subseteq U(N)$ for some $N\ge 2$ and consider any infinite volume limit of four-dimensional lattice gauge theory with gauge group $G$ at inverse coupling strength $\beta$. Let $\gamma_{R,T}$ denote a rectangular loop of length $T$ and breadth $R$. Prove that as $\beta \to \infty$, there are sequences $\ep=\ep(\beta)\to 0$ and $c = c(\beta)\to \infty$, and a nonzero constant $d$, such that for any $R$ and $T$,  
\begin{equation}\label{wlimit}
\log \smallavg{W_{\gamma_{R/\ep, T/\ep}}} = -c(R+T) - dRT + o(1).
\end{equation}
\end{problem}
Note that $R+T$ is the limiting perimeter of the rectangular loops after scaling by $\ep$, and $RT$ is the limiting area. 
Since $c\to \infty$,  the above conjecture says that Wilson's area law  holds in the continuum only after subtracting off (`renormalizing away') the first term when taking the limit. That is, the logarithm of the Wilson loop expectation $\smallavg{W_{\gamma_{R,T}}}$ in the continuum should be defined as
\[
\log \smallavg{W_{\gamma_{R,T}}} := \lim_{\beta\to\infty}(\log \smallavg{W_{\gamma_{R/\ep, T/\ep}}} + c(R+T)). 
\]
With this definition, the string tension of the continuum theory (as defined in Section \ref{confinesec}) is the number $d$ in \eqref{wlimit}. 

Since none of the above has been proved, it is not clear to me whether the renormalization term $c(R+T)$ is indeed necessary. It seems entirely possible that the limit in \eqref{wlimit} holds without the renormalization term. 


The next section gives a brief summary of existing rigorous results on the problem of constructing continuum limits of lattice gauge  theories in various dimensions. 

\section{Review of the mathematical literature}\label{litsec}
There is a long and quite successful development of two-dimensional Yang--Mills theories in the mathematical literature. The two-dimensional Higgs model, which is $U(1)$ Yang--Mills theory with an additional Higgs field, was constructed by \citet{bfs79, bfs80, bfs81} and further developed by \citet{bs83}. Building on an idea of \citet{bralic80}, \citet{gks89} formulated a rigorous mathematical approach to performing calculations in two-dimensional Yang--Mills theories via stochastic calculus. Different ideas leading to the same goal were implemented by \citet{driver89a, driver89b} and \citet{kk87}. The papers of \citet{driver89a, driver89b} made precise the idea of using objects called lassos to define the continuum limit of Yang--Mills theories. Explicit formulas for Yang--Mills theories on compact surfaces were obtained by \citet{fine90, fine91} and \citet{witten91, witten92}. All of these results were generalized and unified by \citet{sengupta92, sengupta93, sengupta97} using a stochastic calculus approach. 

Yet another approach was introduced by \citet{levy03, levy10}, who constructed   two-dimensional Yang--Mills theories as random holonomy fields. A random holonomy field is  a stochastic process indexed by curves on a surface, subject to boundary conditions, and behaving under surgery as dictated by a Markov property. L\'evy's framework allows parallel transport along more general curves than the ones considered previously, and makes interesting connections to topological quantum field theory. A relatively non-technical description of this body of work is given in~\citet{levy11}. 

Recently, \citet{nguyen15} has established the mathematical validity of the perturbative approach to 2D Yang--Mills theory by comparing its predictions with rigorous results obtained by the approaches outlined above. Another important body of recent work consists of the papers of \citet{levy11b}, \citet{driveretal}, \citet{driveretal2} and \citet{driver17}, who establish the validity of the Makeenko--Migdal equations for 2D Yang--Mills theories \cite{makeenkomigdal79} by a number of different approaches.

Euclidean Yang--Mills theories in  dimensions three and four have proved to be more challenging to construct mathematically. At sufficiently strong coupling (that is, small $\beta$), a number of conjectures about lattice gauge theories in arbitrary dimensions --- such as quark confinement and the existence of a positive self-adjoint transfer matrix --- were rigorously proved by \citet{os78}. An expansion of partition functions of lattice gauge theories as asymptotic series in the dimension of the gauge group was proposed by \citet{thooft74}, leading to a large body of work. The papers of~\citet{chatterjee15}, \citet{cj16}, \citet{bg16} and \citet{jafarov16} contain some recent advances on 't Hooft type expansions and connections with gauge-string duality at strong coupling. Confinement  and deconfinement in three- and four-dimensional lattice gauge theories at weak coupling were investigated by \citet{guth80}, \citet{frohlichspencer82},  \citet{gopfertmack82} and \citet{borgs88}. 

None of the above techniques, however, help in constructing the continuum limit. The problems posed in Section \ref{contsec} have not been mathematically  tractable. An alternative route is the method of phase cell renormalization. In this approach, one starts with a lattice gauge theory on $\ep\zz^n$ for some small $\ep$. Choosing an integer $L$, the theory is then `renormalized' to yield an `effective field theory' on the coarser lattice $L\ep\zz^n$, which is just another lattice gauge theory but with a more complicated action. A survey of the various ways of carrying out this renormalization step is given in \cite[Chapter 22]{gj87}. The process is iterated to produce effective field theories on $L^k\ep\zz^n$ for $k=2,3,\ldots$, until the lattice spacing $L^k\ep$ attains macroscopic size (for example, becomes greater than $1$). Note that the macroscopic effective field theory obtained in this way is dependent on $\ep$. The goal of  phase cell renormalization is to show that the effective field theory at the final macroscopic scale converges to a limit as $\ep\to 0$. Usually, convergence is hard to prove, so one settles for subsequential convergence by a compactness argument. The existence of a convergent subsequence is known as ultraviolet stability, for the following reason. The approximation of an Euclidean Yang--Mills theory by a lattice gauge theory on a lattice with spacing $\ep$ is analogous to truncating a Fourier series at a finite frequency, which grows as $\ep\to 0$. In this sense, a lattice approximation is an ultraviolet cutoff (ultraviolet = high frequency), and  the compactness of the effective field theories is ultraviolet stability, that is, stability with respect to the cutoff frequency.

A notable success story of phase cell renormalization is the work of \citet{king86a, king86b}, who established the existence of the continuum limit of the three-dimensional Higgs model. The continuum limit of pure $U(1)$ Yang--Mills theory (that is, without the Higgs field) was established earlier by \citet{gross83}, but with a different notion of convergence. Gross's approach was later used by \citet{driver87} to construct a continuum limit of 4D $U(1)$ lattice gauge theory.

Ultraviolet stability of three- and four-dimensional non-Abelian lattice gauge theories by phase cell renormalization, as outlined above, was famously established by \citet{balaban83, balaban84a, balaban84b, balaban84c, balaban84d, balaban85a, balaban85b, balaban85c, balaban85d, balaban85e, balaban87, balaban88, balaban89a, balaban89b} in a long series of papers spanning six years.  A somewhat different approach, again using phase cell renormalization, was pursued by \citet{federbush86, federbush87a, federbush87b, federbush88, federbush90} and \citet{fw87}. 


Phase cell renormalization is not the only approach to constructing Euclidean Yang--Mills theories. \citet{mrs93} formulated a program of directly constructing Yang--Mills theories in the continuum instead of using lattice theories. The main idea in \cite{mrs93} was to regularize the continuum theory by introducing a quadratic term in the Hamiltonian. The problem with this regularization is that it breaks gauge invariance. The problem is taken care of by showing that it is possible to remove the quadratic term and restore gauge invariance by taking a certain kind of limit of the regularized theories.

In spite of the remarkable achievements surveyed above, there is yet no construction of a continuum limit of a lattice gauge theory in any dimension higher than two where Wilson loop variables have been  shown to have nontrivial behavior. The standard approach of regularizing Wilson loop variables by phase cell renormalization has not yielded definitive results. In a recent series of papers, \citet{cg13, cg15, cg17} and \citet{gross16, gross17} have proposed a method of regularizing connection forms in $\rr^3$ by letting them flow for a small amount of time according to the Yang--Mills heat flow (see also the papers of \citet{luscher10, luscher10b} for a similar idea). It will be interesting to see whether this new approach, in combination with some ideas developed in the paper~\cite{chatterjee16}, can  lead to the construction of nontrivial three-dimensional Euclidean Yang--Mills theories with non-Abelian gauge groups. 


\section*{Acknowledgments}
I thank Erik Bates, David Brydges,  Persi Diaconis, Len Gross, Jafar Jafarov, Erhard Seiler, Scott Sheffield, Steve Shenker, Tom Spencer and Edward Witten   for many valuable conversations and comments.

\end{document}